\documentclass[a4paper,11pt]{amsart}
\usepackage{hyperref}
\usepackage[usenames,dvipsnames,svgnames,table]{xcolor}
\usepackage[all]{xy}
\usepackage{enumitem}
\usepackage[normalem]{ulem}
\newcommand{\tsk}[1]{\textcolor{YellowOrange}}

\usepackage{tikz}
\usepackage{tikz-cd}
\usepackage{graphicx}
\usepackage{pgf,tikz}
\usetikzlibrary{arrows}
\usepackage[left=2cm,right=2cm, top=3cm, bottom=3cm]{geometry}

\makeatletter
\def\@endtheorem{\endtrivlist}% NEW
\makeatother
\usepackage[active]{srcltx}
\usepackage{amsmath}
\usepackage{amssymb}
\usepackage{amscd}
\usepackage{amsthm}
\usepackage{mathrsfs}  

\usepackage{nicefrac}

\newtheorem{teo}{Theorem}[section]
\newtheorem{defin}[teo]{Definition}
\newtheorem{prop}[teo]{Proposition}
\newtheorem{cor}[teo]{Corollary}

\theoremstyle{definition}

\newtheorem{remark}[teo]{Remark}

\newtheorem*{conj*}{Conjecture} %l0ambiente è congettura ma senza il numero 

\newtheoremstyle{dico}% name of the style to be used
 {\baselineskip}   % ABOVESPACE
  {\topsep}   % BELOWSPACE
  {}  % BODYFONT
  {0pt}       % INDENT (empty value is the same as 0pt)
  {} % HEADFONT
  {.}         % HEADPUNCT
  {5pt plus 1pt minus 1pt} % HEADSPACE
  {}          % CUSTOM-HEAD-SPEC
\theoremstyle{dico}

\numberwithin{equation}{section}

\newcommand{\ra}{\rightarrow}
\newcommand{\C}{\mathbb{C}}
\newcommand{\R}{\mathbb{R}}
\newcommand{\Zeta}{{\mathbb{Z}}}

\newcommand{\QQ}{{\mathbb{Q}}}
\newcommand{\meno}{^{-1}}

\newcommand{\vacuo}{\emptyset}

\newcommand{\Aut}{\operatorname{Aut}}

\newcommand{\End}{\operatorname{End}}

\newcommand{\om}{\omega}

\renewcommand{\phi}{\varphi}

\newcommand{\PP}{\mathbb{P}^1}

\newcommand{\OO}{\mathcal{O}}

\newcommand{\sieg}{\mathfrak{S}}

\newcommand{\Z}{\mathsf{Z}}
\newcommand{\perdue}{^{\otimes 2}}

\newcommand{\sym}{{\operatorname{Sym}}}

\newcommand{\ag}{\mathcal{A}_g}
\newcommand{\mg}{\mathcal{M}_g}

\newcommand{\mihi}[1]{}

\newcommand{\Nm}{\operatorname{Nm}}
\newcommand{\Pic}{\operatorname{Pic}}

\begin{document}
 
\pagestyle{myheadings}

\title{Shimura subvarieties via endomorphisms}

\author{Irene Spelta}

\address{Universit\`a degli Studi di Pavia, Dipartimento di Matematica,  Via Ferrata 5  \\ 27100 Pavia, Italy  \noindent}
\email{irene.spelta01@universitadipavia.it}

\thanks{
The author was partially supported by MIUR PRIN 2017
``Moduli spaces and Lie Theory'',  by MIUR ``Programma Dipartimenti di Eccellenza
(2018-2022) - Dipartimento di Matematica F. Casorati
Universit\`a degli Studi di Pavia '' and by INdAM (GNSAGA)}
\date{}

	\begin{abstract}
	We show the existence of two new Shimura subvarieties of $\mathcal{A}_2, \mathcal{A}_3$ generically contained in the Torelli locus. They provide the first examples of Shimura subvarieties obtained by means of Jacobians carrying non-trivial endomorphisms not directly induced by the automorphisms of the curves. We also obtain a new example of a Shimura subvariety of $\mathcal{A}_4$ generically contained in the Prym locus.  
	\end{abstract}
	\maketitle
	\setcounter{section}{-1}

\section{Introduction}
	Let us denote by $\mg$ the coarse moduli space of smooth projective genus $g$ curves and by $\ag$ the one of principally polarized $g$-dimensional abelian varieties. The Torelli map $t: \mg\ra \ag$ sends the class $[C]\in\mg$ to the class of the Jacobian variety of $C$ together with the polarization induced by the cup product. By Torelli theorem the map $t$ is injective. This justifies the interest on the Zarisky closure of the image $t(\mg)$, usually referred to as \textit{Torelli locus}. 
	
	The study of the Torelli locus has been addressed from different perspective in the literature: one among them concerns Shimura subvarieties $S\subseteq\ag$ which are generically contained in it, i.e. $S\subseteq \overline{t(\mg)}$ such that $S\cap t(\mg)\neq\vacuo$. Shimura varieties are by definition Hodge loci for the natural variation of Hodge structure on $\ag$. They have been characterized by Mumford and Moonen as totally geodesic varieties that carry CM-points (\cite{moonen-linearity-1}, \cite{mumford-Shimura}).
	
	The moduli spaces $\mg$ and $\ag$ are complex orbifolds and the Torelli map is an immersion outside the hyperelliptic locus (see \cite{oort steen}). $\ag$ is the quotient of the Siegel space $\sieg_g$, which is a hermitian symmetric domain, by the action of the symplectic group $\text{Sp}(2g,\R)$. The symmetric metric on $\sieg_g$ induces an orbifold metric on $\ag$. Usually, it is called \textit{Siegel metric}.
	
	 As the genus $g$ grows, one expects the Torelli embedding to be more and more curved with respect to this locally symmetric ambient geometry. In terms of Shimura varieties this expectation is formulated as: 
	\begin{conj*}[Coleman - Oort]\label{congettura}
		For $g\gg0$ there should not exist positive dimensional Shimura subvarieties of $\ag$ generically contained in the Torelli locus. 
	\end{conj*}
Up to genus 7, there are examples of Shimura varieties generically contained in the Torelli locus. They have been found considering families of Galois covers $C\ra C':=C/G$ with fixed monodromy, genera and number of ramification points. Let $\Z$ be the closure in $\ag$ of the locus described by $JC$ for $C$ varying in the family. If \begin{equation}\label{star}
	\dim \sym^2 H^0(C, \om_C)^G= \dim \Z \tag{$ \ast $}
\end{equation}
then $\Z$ is Shimura (see \cite{fgp} for covers of $\PP$ and \cite{fpp} for covers of elliptic curves). Using \eqref{star}, there have been discovered 30+6 counterexamples to the Conjecture. In \cite{cgp} it is shown that, up to genus 100, these are the only examples of Galois covers of $\PP$ satisfying \eqref{star}. We remark that the sufficient  criterion \eqref{star} is proved to be necessary in the case of cyclic Galois coverings of $\PP$ by Moonen in \cite{moonen special} and by Mohajer and Zuo in case of 1-dimensional families of abelian covers of $\PP$ (\cite{mohajer zuo }). See also \cite{fredi} for results in this direction concerning higher-dimensional families of abelian covers. 

Actually, more examples are known in the cases of genus $2,3,4$. Indeed infinitely many new examples of totally geodesic subvarieties (countably many of which are Shimura) were found in \cite{fgs}. There the authors considered two fibrations involving the 6 examples of families of Galois coverings of elliptic curves yielding Shimura varieties described in \cite{fpp}. One of these fibrations is given by the Prym map. In fact, the paper \cite{fgs} generalizes a result of Grushevsky and M\"oller (\cite{Grushevsky moller }) which studies the fibres of the Prym map of one among the 6 examples of \cite{fpp} to show that they are totally geodesic. 
	
In \cite{cfgp} the authors reformulated the Coleman-Oort conjecture for the Prym loci, namely they asked if there exist positive dimensional Shimura subvarieties generically contained in the Prym loci $\mathcal{P}(\mathcal{R}_{g,0})$ and $\mathcal{P}(\mathcal{R}_{g,2})$. As usual, we denote by $\mathcal{R}_{g,r}$ the moduli space parametrizing double covers $\pi: C\ra C'$ of genus $g$ curves $C'$ ramified in $r\geq0$ points and $\mathcal{P}$ is the corresponding Prym map. $\mathcal{P}$ sends $[\pi]$ to the isomorphism class of its associated Prym variety. Analogous investigation  has then been done in \cite{fredi gross} and \cite{fredi gross moj}  on  $\mathcal{P}(\mathcal{R}_{g,r})$ for all $r$ and, recently, also in case of covers $\pi$ of arbitrary degree \cite{gross moh}. This kind of problems naturally arises since the second fundamental form of the Prym map has very similar behaviour to the one of the Torelli map $t$. Many examples of Shimura varieties generically contained in the different Prym loci were found studying families of Galois covers of $\PP$ and requiring that the codifferential of the involved Prym map be an isomorphism at the generic point of the family (condition B of \cite{cfgp},  \cite{fredi gross}, \cite{fredi gross moj}). Notice that, as for \eqref{star}, B is only a sufficient condition. 

We stress that all the examples of Shimura varieties, both in the case of the Torelli locus or in the ones of Prym loci, have been constructed by considering families of Galois covers of $\PP$ or of elliptic curves. Thus the automorphism group of the curves, hence of the Jacobians, plays a crucial role. 

Very little is known about the existence of families of curves whose Jacobians are acted on by a large ring of endomorphisms not induced by automorphisms of the curves. The first examples are due to  Tautz-Top-Verberkmoes, Mestre and Brumer. Then Ellenberg presented a general procedure to construct curves whose Jacobians have large endomorphism algebras, getting almost all the examples already known in the literature at his moment (see \cite{ellemberg} and the references therein). In particular for every odd prime number $p\geq 7$ he found 3-dimensional families of curves of genus $(p-1)/2$ whose endomorphism algebra contains the real cyclotomic field $\QQ(\xi_p+\xi_p\meno)$, with $\xi_p^p=1$.

This is our starting point. Indeed in these note, we show the existence of two new Shimura subvarieties of $\mathcal{A}_2$, resp. $\mathcal{A}_3$, of dimension 2, resp. 3, which are generically contained in the Torelli locus. These are the first known examples of Shimura varieties which involve Jacobians whose endomorphism algebra is non-trivial and the endomorphisms are not induced by automorphisms of the curves.

Our construction is inspired by the families already studied by Albano and Pirola in \cite{albano pirola} to disprove a conjecture of Xiao relating the relative irregularity and the genus of the general fiber of a fibration. Indeed in \cite{albano pirola} the authors consider three fibrations whose fibres are certain cyclic \'etale covers of hyperelliptic genus $g$ curves $\pi:C\ra H:=C/\langle\sigma\rangle$ of prime odd degree $p$ and they show that they violate Xiao's conjecture.

It turns out that the curves $C$ occurring in such families are dihedral covers of $\PP$ and that the lift of the hyperelliptic involution of $H$ to $C$ gives an intermediate quotient $C\xrightarrow{2:1} D$. The curves $D$ map $p:1$ to $\PP$ but the covering is non-Galois. These curves are special in moduli. Indeed the automorphism $\sigma$ induces a non-trivial automorphism of $JD$ (not preserving the polarization on $JD$). It follows that $\End_0(JD)$ contains at least $\QQ(\xi_p+\xi_p\meno)$. For details on such covers, we refer to \cite{ries}. 

This construction gives a map \begin{equation*}
\psi: \mathcal{RH}_{g}(p)\ra \mathcal{M}_{\frac{1}{2}(g-1)(p-1)}
\end{equation*}
from the moduli space $\mathcal{RH}_{g}(p)$ parametrizing covers $\pi:C\ra H$, to the moduli space of curves $D$ obtained from $\pi$ as above. By Riemann-Hurwitz formula, their genus is equal to $\frac{1}{2}(g-1)(p-1)$. 

Our results are the following:
\begin{teo}\label{teorema 1}
	For $g=2$ and $p=5$, resp. $p=7$, the image through $t$ of the image of $\psi$ yields a new Shimura subvariety of $\mathcal{A}_2$, resp. $\mathcal{A}_3$, of dimension 2, resp. 3, generically contained in the Torelli locus. 
\end{teo}
In this way, we provide the first two examples of Shimura subvarieties of $\mathcal{A}_2, \mathcal{A}_3$ generically contained in the Torelli locus without using the condition \eqref{star}.

The families analysed in Theorem \ref{teorema 1} can also be studied by looking at their image in the Prym locus. Let us consider the diagram
\begin{equation*}
\begin{tikzcd}
C\arrow[r,"p:1"]\arrow[rr,bend left]&H \arrow[r,"2:1"]&\PP=C/D_p.
\end{tikzcd}
\end{equation*}
  It is known by \cite{ries} that the Prym variety $P(C, H)$ is isomorphic to $JD\times JD$. We have the following consequence:
\begin{prop}
	The Prym varieties of the two families of Theorem \ref{teorema 1} yield Shimura subvarieties of $\mathcal{A}_4$, resp. $\mathcal{A}_6$, of dimension 2, resp. 3, generically contained in the Prym locus.
\end{prop}
Finally, we show that for the one with data $g=2$ and $p=5$ the codifferential of the Prym map is not an isomorphism. Therefore it is not one of the examples found in \cite{gross moh}. In particular, we have the following
\begin{cor}
	The condition B is only sufficient for a family of Galois covers of $\PP$ to yield a Shimura subvariety generically contained in the Prym locus. 
\end{cor}

{\bfseries \noindent{Acknowledgements}}. The author thanks G.P. Pirola for his inspiring question, A. Ghigi for the right reference at the right time and P. Frediani for her very careful reading and her valuable comments. 

	\section{Hyperelliptic coverings and Endomorphism of Jacobians}\label{sezione 1}

	In this section, we overview known results concerning cyclic coverings of hyperelliptic curves and their associated Prym varieties. Moreover, we describe the endomorphism algebras of the Jacobians of the curves occurring through this construction.
	
	 Let $H$ be a hyperelliptic curve of genus $g$ and $\pi:C\ra H$ an unramified cyclic cover of odd prime order $p$. By general theory on cyclic \'etale coverings, $\pi$ corresponds to the pair $(H, L )$, where $L\in \Pic^0(H)$ satisfies $L^p=\OO_H$. Let $\mathcal{RH}_{g}(p)$ be the coarse moduli space of pairs $(H, L )$. We will call such elements \textit{hyperelliptic coverings }. By Riemann-Hurwitz formula we get $g(C)=g_C=p(g-1)+1. $ Let us denote by $\sigma$ a fixed generator of the cyclic group $\mathbb{Z}/p$ acting on $C$. We borrow from \cite{ries} what follows.  
	\begin{prop}\label{gruppo diedrale}
		Let $(H, L)$ be an element of $\mathcal{RH}_{g}(p)$. The hyperelliptic involution of $H$ lifts to an involution $\tau $ on $C$ and the covering $C\ra \PP$ is Galois with Galois group $G=D_p$ generated by $\sigma$ and $\tau$ (here $D_p$ stands for the dihedral group of order $2p$). 
	\end{prop}
Let $C\ra D$ be the quotient of $C$ by the lift of the hyperelliptic involution $\tau$. We look at the following diagram: \begin{equation}\label{diagramma famiglie}
\begin{tikzcd}
C\arrow{d}{2:1}\arrow{r}{p:1}\arrow{dr}&H\arrow{d}{2:1}\\
D\arrow{r}{p:1}&\PP\cong C/D_p.
\end{tikzcd}
\end{equation} 

The cover $D\xrightarrow{p:1}\PP$ is non-Galois: over every Weiestrass point for $H$ the map $D\xrightarrow{p:1}\PP$ has $1+(p-1)/2$ points. One of this point is \'etale while the others have order 2. The Galois closure of such a map is the cover $C\ra C/D_p\cong \PP.$  By Riemann-Hurwitz formula, we have  $g(D)=g_D=(g-1)(p-1)/2. $ 

It is well known that a pair $(H,L)$ determines an abelian variety $P(H, L)$ defined as the connected component containing the origin of the Norm map $\Nm: JC\ra JH$. It carries a polarization $\Xi$ induced by restriction of that of $JC$ and it is usually referred to as the \textit{Prym variety} of the covering $C\ra H$. In the case of hyperelliptic \'etale covers we have the following:
\begin{teo}[Ries \cite{ries}]\label{iso Ries}
	There is an isomorphism of polarized abelian varieties 
	\begin{equation}\label{decomposizione Prym}
		P(H, L)\cong JD\times JD.
	\end{equation}
Hence  $JC$ is isogenous to $ JH\times JD\times JD.$ Furthermore the endomorphism $\sigma^*+(\sigma\meno)^*$ of $JC$ induces a non-trivial endomorphism of $JD$ for $p>3$. In particular, the endomorphism $\sigma^{\frac{p+1}{2}}+(\sigma\meno)^{\frac{p+1}{2}}$ restricts to an automorphism of $JD$ which does not preserve the polarization on $JD$. 
\end{teo}

\begin{defin}
	Let $F$ be a totally real number field, i.e. $F$ is generated over $\QQ$ by a root of an integer polynomial, whose roots are all real. A curve $X$ has real multiplication by $F$ if $g(X)=[F:\QQ]$ and  if $F\hookrightarrow \End_0(JX)$. 
\end{defin}
\begin{cor}\label{endomorfismi}
	When $g=(p-1)/2$ and $p>3$ the curve $D$ has real multiplication by $\QQ(\xi_p+\xi_p\meno)$, where $\xi_p^p=1, \xi_p\neq 1$. 
\end{cor}

Using Proposition \eqref{gruppo diedrale}, we can define the morphism \begin{equation*}
	\psi: \mathcal{RH}_{g}(p)\ra \mathcal{M}_{\frac{1}{2}(g-1)(p-1)}
\end{equation*}
which sends the isomorphism class of a covering $C\ra H$ to the isomorphism class of the associated quotient curve $D$. This is well-defined since any two lifts of the hyperelliptic involution are conjugated in $\Aut(C)$.

On the other hand we can consider the Prym map 
\begin{equation}\label{mappa di Prym}
	\mathcal{P}: \mathcal{RH}_{g}(p)\ra \mathcal{A}_{(g-1)(p-1)} , \qquad (H,L)\mapsto (P(H,L), \Xi)
\end{equation}
which, thanks to Theorem \ref{iso Ries},  can be thought as $T\circ \psi$, where $T$ is defined as follows  \begin{align*}
	T: \mathcal{M}_{\frac{1}{2}(g-1)(p-1)}&\ra \mathcal{A}_{(g-1)(p-1)}\\
	D&\mapsto JD\times JD.
\end{align*}

 By Torelli Theorem, we can observe that the fibres of $\mathcal{P}$ have the same dimension of those of $\psi$. They have been described by the following 
\begin{prop}[Albano-Pirola, \cite{albano pirola}]\label{fibre psi}
	The map $\psi$ has finite fibres if and only if $p\geq 7$ or $p=5$ and $g\geq 3$ or $p=3$ and $g\geq 5$.
\end{prop}
In \cite{albano pirola}, the authors used Proposition \ref{fibre psi} to look for positive dimensional fibres for $\psi$. In those cases, an irreducible component of the fiber $\psi\meno (D) $ gives families of curves $C$ (as in diagram \eqref{diagramma famiglie}) whose decomposition of the Jacobians contains the fixed part $JD\times JD$. Their goal was to study the relative irregularity of the fibrations obtained in this way in order to compare it with the genus of the general fiber. Indeed they showed that:
\begin{teo}[\cite{albano pirola}]\label{fibrations AP}
	There exist three fibrations which disprove Xiao's conjecture on relative irregularity. They have data $g=2 $ and $p=5$, $g=4$ and $p=3$ and $g=3$ and $p=3$.
\end{teo} 
\begin{remark}
 Actually, the first Xiao fibration, namely the first fibration which violates Xiao's conjecture, was constructed by Pirola in \cite{pirola}. It consists of a 3-dimensional family of curves which are Galois triple covering of elliptic curves. This family is very interesting. Indeed in \cite{fpp} it is showed that it yields a Shimura 3fold in $\mathcal{A}_4$ and in \cite{fgs} that, through its Prym map, it is fibred in infinitely many totally geodesic fibres, countably many of which are Shimura. One among these Shimura fibres is the 1-dimensional family (12) of \cite{fgp} given by genus 4 curves which are $\mathbb{Z}/6$-Galois covers of $\PP$. Furthermore, the family of \cite{pirola} has been got in \cite{gross spelta} where the authors study Xiao fibrations considering Galois covers $C\ra H$ under less restrictive assumptions (i.e. any degree and any value for $r$). 
\end{remark}
Although we no longer deal with Xiao's conjecture, the families described in \cite{albano pirola} turn out to be suitable for our purposes. Indeed up to now we have shown that, when $p>3$, the construction of diagram \eqref{diagramma famiglie}, together with Corollary \ref{endomorfismi}, produces families of curves $D$ of genus $(g-1)(p-1)/2$ with $\QQ(\xi_p+\xi_p\meno)\hookrightarrow\End_0(JD)$. 

When the fibres of $\psi$ are finite, these families have dimension $2g-1$, i.e. the dimension of $\mathcal{RH}_{g}(p)$. In particular, when $g=2$ and $p\geq7$ we recover the 3-dimensional families of curves whose Jacobians are acted on by a large ring of endomorphism found in case 1 of the main Theorem of \cite{ellemberg}. 
On the other hand, assuming that the fibres of $\psi$ are positive dimensional,  a straightforward computation shows that we get a unique family worth of interest in the sense of its endomorphism algebra. It is the one already studied by Albano and Pirola with data $g=2$, $p=5$ and it has dimension 2. 

In the next section, we want to show that there are two cases when the Jacobians $JD$ of the curves $D$ occurring in such families yield Shimura subvarieties of $ \mathcal{A}_{\frac{1}{2}(g-1)(p-1)} $ generically contained in the Torelli locus. Moreover, thanks to Theorem \ref{iso Ries}, we obtain interesting consequences on Prym loci too. 

\section{Shimura Subvarieties}
Let $\Lambda$ be a rank $2g$ lattice, $Q: \Lambda\times \Lambda\ra \Zeta$ an alternating form of type $(1,1,...,1)$. Let $K$ be a field $\mathbb{Q}\subseteq K\subseteq\mathbb{C}$ and set $\Lambda_K:=\Lambda\otimes_\Zeta K$. Denote by $A_J$ the quotient $\Lambda_\mathbb{R}/\Lambda$ provided with the complex structure $J$ and the polarization $Q$. It yields an element in $\ag$. As usual, we will interpret $\ag$ as the orbit space of the action of the symplectic group $\text{Sp}(\Lambda, Q)$ on the Siegel space $\sieg_g$. There is a natural variation of rational Hodge structure on $\ag$: it corresponds to the decomposition of $\Lambda_\C$ into $\pm i$ eigenspaces for $J$. 
\begin{defin}
	A Shimura subvariety of $S\subseteq \ag$ is by definition a Hodge locus of the natural variation of Hodge structure on $\ag$ described above. 
\end{defin}
We will focus on Shimura subvariety of \textit{PEL type}. Fix a complex structure $J_0\in\sieg_g$ and set $F:=\End_\mathbb{Q}(A_{J_0})$, where $\End_\mathbb{Q}(A_{J_0}):=\{f\in \End(\QQ^{2g}):J_0f=fJ_0\}$. We have the following:
\begin{defin}\label{definizione PEL}
	The PEL type Shimura subvariety $S(F)$ is defined as the image in $\ag$ of the connected component containing $J_0$ of the set $\{J\in\sieg_g: F\subseteq\End_\mathbb{Q}(A_J) \}$.
\end{defin}

In view of construction \eqref{diagramma famiglie} and of Corollary \ref{endomorfismi}, it is natural to investigate the subvarieties of $\mathcal{A}_{\frac{1}{2}(g-1)(p-1)}$, contained in the Torelli locus, which are described by the Jacobians $JD$. Our aim is to compare them with the PEL type Shimura subvarieties $S(F_p)$ with $F_p= \mathbb{Q}(\xi_p+\xi_p\meno)$, $p\geq5$.

We have the following:
\begin{teo}\label{nuove shimura}
	The family of curves $D$ of \eqref{diagramma famiglie} with $g=2$ and $p=5$, resp. $p=7$, gives rise to a Shimura subvariety of PEL type $S(F_5)\subset \mathcal{A}_2$, resp. $S(F_7)\subset \mathcal{A}_3$, generically contained in the Torelli locus of dimension 2, resp. 3.
	\begin{proof}
		(1). Let us start with $g=2, p=5$ (precisely with one of the families studied in \cite{albano pirola}). The genus 2 curves $D$ have real multiplication by  $ F_5= \mathbb{Q}(\xi_5+\xi_5\meno) $. Indeed $F_5$ is a totally real field with $e:=[F_5:\mathbb{Q}]=2=g_D$. Let $S(F_5)$ be the Shimura subvarieties of $\mathcal{A}_{2}$ of PEL type defined by $F_5$ as in Definition \ref{definizione PEL}. Hence $S(F_5)$ arises as the quotient of $g_D$ copies of the Siegel space $\sieg_{g_D/e}=\sieg_1$, as showed in \cite[Theorem 4.16]{moonen-oort} (see also \cite[subsection 2.3]{jong Noot}). Therefore \[\dim S(F_5)= \dfrac{g_D(g_D+e)}{2e}=2.	\]
		Now let $X$ be the closure of the locus of $JD$ for $D$ varying in the family of Theorem \ref{fibrations AP} with $g=2$ and $p=5$. Then clearly $X$ is an irreducible subvariety of $S(F_5)$ and $\dim X= \dim S(F_5)=2$. Thus $X$ must coincides with the PEL subvariety of $\mathcal{A}_2$ identified by $F_5$. Therefore it is Shimura.
		
		(2). Now take $g=2, p=7$. In this case we have \[F_7= \mathbb{Q}(\xi_7+\xi_7\meno), \ \ e=g_D=3. \] The same computation as in the previous situation gives $\dim S(F_7)=3$. Since $JD$ varies in a 3-dimensional family we get the PEL subvariety of $\mathcal{A}_3$ identified by $F_7$. This concludes the proof. \\
	\end{proof}
\end{teo}
 \begin{remark}
 	This Theorem determines the first two examples of Shimura subvarieties of $\ag$ generically contained in the Torelli locus which are not obtained by means of families of Galois coverings of $\PP$ satisfying the sufficient condition \eqref{star}. Actually, we need to check that our family (1), resp. (2), is different from those obtained in \cite{fgp} and in \cite{fpp}. Genus and dimension tell us that they might coincide only with family (26), resp. family (27), of \cite{fgp}. Family (26) is the bielliptic locus in genus 2 and family (27) is strictly contained in the bielliptic locus in genus 3. The generic curve of these two families carries a completely decomposable Jacobian. Using Mumford's list of non-trivial endomorphism rings of abelian varieties (\cite{mumford libro}) and comparing the endomorphism algebras of families (26), (27) with the ones of our families (1) and (2), we conclude.  	
 	Finally, notice that \cite[Theorem 1.1]{cvdgt} states that families (26) and (1) are the unique irreducible 2-dimensional subvarieties of $\mathcal{A}_2$ whose generic point has non-trivial endomorphisms, namely $\Zeta\subsetneq\End_0(JC)$. 
 	%Furthermore Theorem \ref{iso Ries} states that the curves of family (1), resp. (2), of Theorem \ref{nuove shimura} carry the automorphism $\sigma^3+\sigma^2$, resp. $\sigma^4+\sigma^3$. None of the two are involutions, hence they do not belong to the automorphism groups of families (26) and (27) which are isomorphic to $\Zeta/2\times\Zeta/2.$ Hence we conclude.  QUESTO è ERRATO: L'AUTOMORFISMO DELLA JACOBIANA NON DISCENDE AD UN AUTOMORFISMO DELLA CURVA PERCHè NON PRESERVA LA POLARIZZAZIONE quindi non è questo il modo per confrontarli. 
 \end{remark}
\begin{remark}
		For other values of $g$ and $p$, the image of $\psi$ still produces subvarieties of $t(\mathcal{M}_{g_D})$. In these cases, such subvarieties are strictly contained in the Shimura subvarieties of PEL type arising from the corresponding totally real field $F_p= \mathbb{Q}(\xi_p+\xi_p\meno)$. This agrees with \cite[Corollary 4.12]{moonen-oort} where it is stated that if $g\geq4$ and $S\subseteq\ag$ is a Shimura subvariety of PEL type obtained from a totally real field of degree $g$, then $S$ is not generically contained in the Torelli locus. This result is obtained in \cite{jong Noot} for $g>4$ and in \cite{bainbridge Moller} for $g=4$. Therefore, letting $g$ and $p$ vary, we should not expect the existence of other families as (1) and (2) yielding Shimura varieties of PEL type generically contained in higher genus Torelli locus. Notice that our examples prove that the bound $g\geq 4$ of \cite[Corollary 4.12]{moonen-oort} is sharp.
\end{remark}
Let us now look at the loci described by the Prym varieties of the two families studied in Theorem \ref{nuove shimura}, namely at the families of $P(C,H)$ in $\mathcal{P}(\mathcal{RH}_{g}(p))\subset \mathcal{A}_{(g-1)(p-1)}$. 

Let us interpret an element $C\ra H$ of such families as a Prym datum $(D_p, \theta, \sigma)$ (we borrow the notation of  \cite{cfgp}, \cite{fredi gross}, \cite{fredi gross moj}): $D_p$ is the Galois group acting on $C$ with quotient $\PP$, $\sigma$ gives the intermediate quotient $C\ra H= C/\langle\sigma\rangle $ and $\theta$ is the monodromy map. Moreover setting $V:=H^0(C,\om_C)$ and letting $\sigma$ acts on $V$, we can decompose \[V=V_+\oplus V_{-}\quad \text{where}\quad V_+\cong H^0(C,\om_C)^{\langle\sigma\rangle}= H^0(H, \om_H)\quad \text{and}\quad V_{-}\cong\bigoplus_{i=1}^{p-1} H^0(H, \om_H\otimes L^i ). \]
By \cite[Proposition 4.1]{lange ortega}, the codifferential of the Prym map \eqref{mappa di Prym} at a point $(H, L)\in \mathcal{RH}_{g}(p)$ can be identified with the multiplication map \[m: \sym^2 V_{-}\ra H^0(H, \om_H\perdue).\]
For our purposes, we are interested in the restriction of $m$ to $\sym^2 V_{-}^{\langle\sigma\rangle}$. Indeed $\langle\sigma\rangle \subseteq \Aut(JC)$ acts on $JC$ and thus on $P(C,H)$. Therefore the image of $\mathcal{P}$ is contained in the locus of ${(g-1)(p-1)}$-abelian varieties with $\Zeta/p$ automorphisms. Hence the codifferential of $\mathcal{P}$ is given by the multiplication map \begin{equation}\label{codifferential Prym}
	 m: (\sym^2\bigoplus_{i=1}^{p-1} H^0(H, \om_H\otimes L^i ))^{\langle\sigma\rangle} \ra H^0(H, \om_H\perdue).
\end{equation}
Theorems 3.2 of \cite{cfgp}, \cite{fredi gross}, \cite{fredi gross moj} use the map $m$ in \eqref{codifferential Prym} to give a sufficient criterion for a family of Galois covers $C\xrightarrow{2:1} C/\langle \sigma \rangle \ra \PP$, identified by the Prym datum $(G, \theta, \sigma)$, to yield a Shimura subvarieties of $\mathcal{A}_{g-1+\frac{r}{2}}^\delta$ generically contained in the Prym loci $\mathcal{P}(\mathcal{R}_{g,r})$. Analogous statement is provided in case of higher degree Galois covers $C\ra C/\langle \sigma \rangle$ in \cite[Theorem 4.1]{gross moh}. %(here $G$ denotes the group acting on $C$, $g$ denotes the genus of $H$ and $r\geq0$ the numbers of ramification points for the cover $C\ra H$). 
Indeed the authors show the following:
\begin{teo}[\cite{cfgp}, \cite{fredi gross}, \cite{fredi gross moj}, \cite{gross moh}]
	Let $(G, \theta, \sigma)$ be a Prym datum. If the multiplication map \eqref{codifferential Prym} is an isomorphism at the general point of the family identified by $(G, \theta, \sigma)$, then the closure of the locus of the Prym varieties $P(C,H)$ of the family yield a Shimura subvariety generically contained in the Prym locus.
\end{teo} The assumption of this Theorem is referred to as condition B of \cite{cfgp}, \cite{fredi gross}, \cite{fredi gross moj}, \cite{gross moh}. Hence the authors use this criterion to explicitly construct examples of such varieties.

Obviously, we have that families (1) and (2) yield two Shimura subvarieties of $\mathcal{A}_4$, resp. $\mathcal{A}_6$, generically contained in $\mathcal{P}(\mathcal{RH}_2(5))$, resp. $\mathcal{P}(\mathcal{RH}_2(7))$. This immediately follows from the isomorphism \eqref{decomposizione Prym} and from Theorem  \ref{nuove shimura}. These varieties have dimension 2, 3 respectively.
It thus seems interesting to check if they satisfy or not the aforementioned condition B. Indeed it is only known to be a sufficient criterion. We have the following:
\begin{prop}
	The multiplication map of family (1) is not an isomorphism while the one of family (2) is. 
	\begin{proof}
		The proof is very easy and follows directly from the fact that the fibres of $\mathcal{P}$ have the same dimension of that of $\psi$. In case of family (1), i.e. if $g=2, p=5$, $\psi$ has fibres of dimension 1 (see \cite[Proposition 4.2]{albano pirola}). Therefore \[d\mathcal{P}^*=m: H^0(H, \om_H\otimes L)\otimes H^0(H, \om_H\otimes L^4)\oplus H^0(H, \om_H\otimes L^2)\otimes H^0(H, \om_H\otimes L^3)\ra  H^0(H, \om_H^2) \]
		cannot be surjective. In case of family (2) we have that $d\psi$ is injective and so $d\mathcal{P}$ too (Proposition \ref{fibre psi}). Therefore 
		$m$ is surjective. Furthermore $\dim \sym^2(V_{-})^{\langle\sigma\rangle}=3= \dim H^0(H,\om_H\perdue)$, since $\sym^2(V_{-})^{\langle\sigma\rangle}=H^0(H, \om_H\otimes L)\otimes H^0(H, \om_H\otimes L^6)\oplus H^0(H, \om_H\otimes L^2)\otimes H^0(H, \om_H\otimes L^5)\oplus H^0(H, \om_H\otimes L^3)\otimes H^0(H, \om_H\otimes L^4)$ and $h^0(H, \om_H\otimes L^i)=1, \forall \ i$. Hence we conclude. 
	\end{proof}
\end{prop}
Therefore we can state the following:
\begin{cor}
	The condition B used to produce Shimura subvarieties generically contained in the Prym locus (arising as Galois coverings of the line) is only sufficient. Therefore family (1) yields a new example of Shimura subvarieties of $\mathcal{A}_4$ generically contained in the Prym locus.  
\end{cor}


\begin{thebibliography} {99}
\bibitem{albano pirola} A. Albano and G.P. Pirola, {\em Dihedral monodromy and Xiao fibrations.} Ann. Mat. Pura ed App. 195 (2016), p. 1255-1268.

\bibitem{bainbridge Moller} M. Bainbridge and M. M\"oller,{\em The locus of real multiplication and the Schottky locus.}  J. Reine Angew. Math. 686 (2014), p. 167–186. 

\bibitem{cvdgt} C. Ciliberto, G. van der Geer and M. Teixidor i Bigas, {\em On the number of parameters of curves whose Jacobians possess non-trivial endomorphisms.}
J. Algebraic Geom. 1 (1992), no. 2, p. 215–229. 

\bibitem{cfgp} E. Colombo, P. Frediani, A. Ghigi and M. Penegini, {\em Shimura curves in the Prym locus.} Comm. in Cont. Math. 21 (2019), no. 2, p. 1850009-1850043.

\bibitem{cgp} D. Conti, A. Ghigi and R. Pignatelli, {\em Some evidence for the Coleman-Oort Conjecture.}  arXiv:2102.12349.

\bibitem{ellemberg} J.S. Ellenberg, {\em Endomorphism algebras of Jacobians.} Adv. Math.162 (2001), no. 2, p. 243–271.

\bibitem{fredi} P. Frediani,  {\em Abelian covers and Second Fundamental Form.}  arXiv:2105.07947. 

	\bibitem{fgp} P.~Frediani, A.~Ghigi and M.~Penegini, {\em Shimura
	varieties in the Torelli locus via Galois coverings}. 
	Int. Math. Res. Not. 20 (2015), p. 10595-10623.

	\bibitem{fgs}
	P.~Frediani, A.~Ghigi and I.~Spelta, {\em Infinitely many Shimura varieties in the Jacobian locus for $g \leq 4$}. To appear in Annali SNS. arXiv:1910.13245.

\bibitem{fredi gross} P. Frediani and G. P. Grosselli, {\em Shimura curves in the Prym loci of ramified double covers.}  arXiv:2007.09646. 

\bibitem{fredi gross moj} P. Frediani, G.P. Grosselli and A. Mohajer, {\em Higher dimensional Shimura varieties in the Prym loci of ramified double covers.}  arXiv:2101.09016.

	\bibitem{fpp} P.~Frediani, M.~Penegini and P.~Porru, {\em   Shimura
		varieties in the Torelli locus via Galois coverings of elliptic	curves}. Geometriae Dedicata 181 (2016), p. 177-192.
	
\bibitem{gross moh} G.P. Grosselli and A. Mohajer, {\em Shimura subvarieties in the Prym locus of ramified Galois coverings.} arXiv:2106.05704

\bibitem{gross spelta} G.P. Grosselli and I. Spelta, {\em Explicit analysis of positive dimensional fibres of $\mathcal{P}_{g,r}$ and Xiao conjecture} In preparation.
	
\bibitem{Grushevsky moller } S. Grushevsky and M. M\"oller, {\em Shimura curves within the locus of hyperelliptic Jacobians in genus 3.} Int. Math. Res. Not. IMRN 6 (2016), p. 1603–1639.

\bibitem{jong Noot} A. de Jong and S. Zhang, {\em Generic abelian varieties with real multiplication are
not Jacobians.} Diophantine Geometry, CRM Series 4 (2007), p. 165–172.


\bibitem{lange ortega} H. Lange and A. Ortega, {\em Prym varieties of cyclic coverings.} Geom. Dedicata 150 (2011), p. 391–403.

\bibitem{mohajer zuo } A. Mohajer and K. Zuo, {\em On Shimura subvarieties generated by families of abelian covers of $\PP$}. J. Pure Appl. Algebra 222 (2018), no. 4, p. 931–949.

	\bibitem{moonen-linearity-1} B.~Moonen.  \newblock{\em  Linearity
	properties of {S}himura varieties. {I}.}  J. Algebraic
	Geom. 7 (1998), no. 3, p. 539-567.

 \bibitem{moonen special} B. Moonen, {\em Special subvarieties arising from families of cyclic covers of the projective line}. Doc. Math. 15 (2010), p. 793–819.

	\bibitem{moonen-oort} B.~Moonen and F.~Oort,{\em   The {T}orelli
		locus and special subvarieties}. {H}andbook of
{M}{oduli} 2 (2013), p. 549-594. 

\bibitem{mumford libro} D. Mumford. {\em Abelian Varieties.} Oxford University Press (1974).

		\bibitem{mumford-Shimura} D.~Mumford. 	{\em A note of {S}himura's	paper ``{D}iscontinuous groups and abelian varieties''. }
 Math. Ann. 181 (1969), p. 345-351.
 
\bibitem{oort steen} F. Oort and J. Steenbrink.	{\em  The local Torelli problem for algebraic curves.}Journ\'ees de G\'eometrie Alg\'ebrique d'Angers, Sijthoff \& Noordhoff, Alphenaan den Rijn (1980), p. 157–204. 
	
\bibitem{pirola} G.P. Pirola, {\em On a conjecture of Xiao.} J. Reine Angew. Math. 431 (1992), p. 75–89.
	
\bibitem{ries} J. Ries, {\em The Prym variety for a cyclic unramified cover of a hyperelliptic Riemann surface.} J. Reine Angew. Math. 340 (1983), p. 59-69. 



	


\end{thebibliography}
\end{document}